# Determination of the asymptotic behavior of the number of natural solutions for certain types of diagonal Diophantine equations

VICTOR VOLFSON

ABSTRACT. We obtain asymptotic upper bounds for the number of natural solutions of the following diagonal Diophantine equations (in a hypercube with side - $N$) in the paper:

$$x_1 = x_2^k + ... + x_s^k,$$

$$x_1^k = x_2^k + ... + x_s^k,$$

$$x_1 = \sum_{j=2}^{s} x_j^{k_j},$$

where $k, s, k_j$ are the natural numbers.

## 1. INTRODUCTION

The n – th degree algebraic Diophantine equation of k variables with integer coefficients has the form: $\sum_{i=0}^{n} F_i(x_1,...x_k) = 0$, where $F_i(x_1,...x_k)$ is the i - th power form of k variables, and $F_0$ is an integer. Thue showed that an irreducible equation $F_n(x_1,...x_k) + F_0 = 0$ ($n \geq 3$, for $k = 2$ and the value $F_0 \neq 0$) has only a finite number of integer solutions, and in a particular case it may not have them at all.

Thue method was developed in the works of Siegel [1], who established on its basis the famous theorem on the finiteness of the number of integer points on a genus curve $g > 0$. This method was extended by Schmidt [2] to the case $k > 2$, which allowed him to obtain a multidimensional generalization of Thue's results on the finiteness of the number of integer solutions of the normative Diophantine equation.

_________________________________________________________________

Quantitative estimates have been obtained for integer solutions of a number of classical Diophantine equations. Baker [3] made effective estimates for the Thue equation of any degree with value $k = 2$. A similar effective analysis for the more general Thue-Mahler equation: $f(x, y) = mp_1^{x_1}...p_s^{x_s}$, where $p_1,..., p_s$ are fixed numbers, was carried out in [4].

Another class of Diophantine equations admitting effective analysis is constituted by superelliptic equations: $y^s = f(x)$, where $s \geq 2$, $f$ is an integer degree function with $n \geq 3$. Effective analysis of the equation: $y^2 = f(x)$ was carried out by Baker [5]. This result was significantly strengthened by Sprinjuk [6]. However, as we see, this analysis is done only for the equations of two variables.

Hardy-Littlewood (CM) circular method [7] allows one to make an upper estimate of the number of natural solutions of various algebraic Diophantine equations with integer coefficients of a diagonal form with a large number of variables. The method was significantly enhanced by Vinogradov [8].

Using the CM in [9], upper bounds were obtained for the number of natural solutions of the diagonal Thue equation: $a_1 x_1^n +...+ a_k^n x_k^n + a_0 = 0$ (for the case when all the coefficients are integers of one sign and $a_0$ is an integer of another sign and not null). The CM was used to obtain an estimate of the number of natural solutions of various homogeneous algebraic Diophantine equations of diagonal form with integer coefficients ($a_0 = 0$) in [10]. The estimate of the number of natural solutions of an inhomogeneous algebraic Diophantine equation of diagonal form was studied in [11].

Estimates of the number of integer (natural) solutions of the Thue equation are found in [12] for values $k > 2$, $n \geq 2$, when the coefficients of the Thue equation are of different signs. An inhomogeneous Diophantine equation of diagonal form with a pronounced variable is also considered there.

Estimates of the number of solutions of Diophantine equations in [9], [10], [11], [12] were carried out mainly for diagonal equations that correspond to the canonical equations of surfaces. The transformations of algebraic Diophantine equations to the diagonal form, preserving the asymptotic number of its integer (natural) solutions, were considered in [13], [14].

We obtain asymptotic upper bounds for the number of natural solutions of the following diagonal Diophantine equations in a hypercube with side - $N$ in the paper:

$$x_1 = x_2^k + \ldots + x_s^k,$$

$$x_1^k = x_2^k + \ldots + x_s^k,$$

$$x_1 = \sum_{j=2}^{s} x_j^{k_j},$$

where $k, s, k_j$ are the natural numbers.

## 2. ESTIMATES OF THE NUMBER OF NATURAL SOLUTIONS OF SERTAIN TYPES OF ALGEBRAIC DIOPHANTINE EQUATIONS OF DIAGONAL FORM

Vinogradov [8] received the following formula for determining the number of representations $n$ as the sum of $s$ natural numbers of $k$ degree:

$$R_s^p(n) = \int_0^1 f^s(x) e^{-2\pi x n} dx, \qquad (2.1)$$

where

$$f(x) = \Sigma_{m=1}^{[n^{1/k}]} e^{2\pi x m^k}, \qquad (2.2)$$

and $[A]$ is the integer value of the number $A$.

Having in mind (2.1), (2.2), the following estimate holds:

$$|R_s^p(n)| = |\int_0^1 f^s(x) e^{-2\pi x n} dx| \leq \int_0^1 |f(x)|^s \, dx, \qquad (2.3)$$

as $|e^{-2\pi x n}| = 1$.

Thus, it is necessary an upper bound of the integral: $\int_0^1 |f(x)|^s \, dx$.

Hua lemma is known [7]. Let be $1 \leq j \leq k$. Then:

$$\int_0^1 |f(x)|^{2^j} \, dx \ll N^{2^j - j + \epsilon}, \qquad (2.4)$$

where value $N = [n^{1/k}]$ (the integer part), and $\epsilon$ is a small real positive number. The $\ll$ icon was introduced by Vinogradov and means that $\int_0^1 |f(x)|^{2^j} \, dx \leq CN^{2^j - j + \epsilon}$, where $C$ is a constant.

Based on the Hua lemma (2.4), we obtain for values $s = 2$ and $j = 1 \leq k$:

$$|R_2^p(n)| \leq \int_0^1 |f(x)|^2 dx << N^{1+\xi}. \tag{2.5}$$

Let's substitute the value $N = [n^{1/k}]$ in (2.5), and then we get:

$$|R_2^p(n)| << n^{1/k+\xi} \tag{2.6}$$

at $k \geq 1$.

Based on the Hua lemma (2.4), we obtain when values $s = 4$ and $j = 2 \leq k$:

$$|R_4^p(n)| \leq \int_0^1 |f(x)|^4 dx << N^{2+\xi}. \tag{2.7}$$

Let's substitute the value $N = [n^{1/k}]$ in (2.7), then we get:

$$R_4^p(n) << n^{2/k+\xi} \tag{2.8}$$

at $k \geq 2$.

Based on the inequality of Cauchy Bunyakovsky (Schwartz) we get for the value $s = 3$:

$$|R_3^p(N)| \leq \int_0^1 |f(x)|^3 dx = \int_0^1 |f(x)|^2 |f(x)| dx \leq (\int_0^1 |f(x)|^4 dx)^{1/2} (\int_0^1 |f(x)|^2 dx)^{1/2}. \tag{2.9}$$

Having in mind (2.4), (2.9) and Hua lemma we obtain:

$$|R_3^p(N)| \leq (\int_0^1 |f(x)|^4 dx)^{1/2} (\int_0^1 |f(x)|^2 dx)^{1/2} << (N^{2+\xi_1})^{1/2} (N^{1+\xi_2})^{1/2} = N^{3/2+\xi}. \tag{2.10}$$

Let's substitute the value $N = [n^{1/k}]$ in (2.10), and then we get:

$$R_3^p(n) << n^{3/2k+\xi} \tag{2.11}$$

at $k \geq 2$.

Let's look at the Diophantine equation:

$$x_1^k + \ldots + x_s^k = n, \tag{2.12}$$

where $s, k, n$ are the natural numbers.

Equation (2.12) does not always have solutions in natural numbers; therefore, the lower bound on the number of these solutions is trivial (equal 0). We will be interested in an asymptotic upper bound for the number of solutions of equation (2.12) in the case when this equation has the indicated solutions.

The number of representations of a natural number $n$, as the sum of $s$ natural numbers of $k$ degree, is related to the number of natural solutions of equation (2.12) as follows:

$$R_s^+(n) = s! R_s^p(n). \tag{2.13}$$

Having in mind (2.13), the asymptotic upper bounds for the number of natural numbers of solutions obtained in formulas (2.6), (2.8) and (2.11) are suitable for equation (2.12) for the corresponding cases.

Based on [7], the asymptotic estimate of the number of natural solutions of equation (2.12) for the value $s > 2^k$ is:

$$|R_s^+(n)| \ll n^{s/k-1}. \tag{2.14}$$

Now we make an asymptotic estimate of the number of natural solutions in a hypercube with side equal $N$ for an equation with an explicit variable:

$$x_1 = x_2^k + \ldots + x_s^k. \tag{2.15}$$

We use the Hardy-Littlewood method for this.

Assertion 1

An asymptotic estimate of the number of natural solutions for an equation with a pronounced variable $x_1 = x_2^k + \ldots + x_s^k$ in a hypercube with a side equal to $N$ (for values $s - 1 > 2^k$) is:

$$|R_s^+(N)| \ll N^{(s-1)/k}.$$

Proof

First, let's look at the equation:

$$x_2^k + \ldots + x_s^k = n. \tag{2.16}$$

Asymptotic estimate of the number of natural solutions of equation (2.16), obtained on the basis of (2.14) for values $s-1 > 2^k$, is equal to:

$$|R_{s-1}^{+}(n)| << n^{(s-1)/k-1}. \qquad (2.17)$$

Based on (2.17), the asymptotic estimate of the number of natural solutions of equation (2.15) in the hypercube with the side $N$, for the values $s-1 > 2^k$, is

$$R_s^{+}(N) = \sum_{n=1}^{N} R_{s-1}^{+}(n) \leq C \sum_{n=1}^{N} n^{(s-1)/k-1}. \qquad (2.18)$$

Having in mind that $n^{(s-1)/k-1}$ is non-decreasing function of $n$ and based on (2.18) we obtain:

$$\sum_{n=1}^{N} n^{(s-1)/k-1} = \int_{1}^{N} t^{(s-1)/k-1} dt + O(N^{(s-1)/k-1}) << N^{(s-1)/k}. \qquad (2.19)$$

Based on (2.18) and (2.19), the asymptotic estimate of the number of natural solutions of equation (2.15) in a hypercube with the side $N$, for values $s-1 > 2^k$, is:

$$|R_s^{+}(N)| << N^{(s-1)/k}, \qquad (2.20)$$

which corresponds to assertion 1

Formula (2.20) is true for values $s-1 > 2^k$. Now we examine the cases when $s-1 \leq 2^k$.

Assertion 2

The asymptotic estimate of the number of natural solutions of the equation $x_1 = x_2^k + x_3^k$ in a hypercube with a side - $N$ is:

$$R_3^{+}(N) << N^{1+1/k+\xi}.$$

The asymptotic estimate of the number of natural solutions of the equation $x_1 = x_2^k + x_3^k + x_4^k + x_5^k$ in a hypercube with a side - $N$ is:

$$R_3^{+}(N) << N^{1+2/k+\xi}.$$

The asymptotic estimate of the number of natural solutions of the equation $x_1 = x_2^k + x_3^k + x_4^k$ in a hypercube with a side - $N$ is:

$$R_4^+(N) \ll N^{1+3/2k+\xi}.$$

Proof

Let's look at the equation:

$$x_1 = x_2^k + x_3^k.$$

Having in mind (2.6), the equation $x_2^k + x_3^k = n$ has the following estimate for the number of natural solutions:

$$|R_2^+(n)| \ll n^{1/k+\xi}.$$

Therefore, we obtain:

$$R_3^+(N) = \sum_{n=1}^{N} R_2^+(n) \le C \sum_{n=1}^{N} n^{1/k+\xi},$$

where C is the constant.

Having in mind that $n^{1/k+\xi}$ is one of a non-decreasing function, we obtain the following estimate of the number of natural solutions in a hypercube with a side - $N$:

$$R_3^+(N) \le C \sum_{n=1}^{N} n^{1/k+\xi} = C \int_{1}^{N} t^{1/k+\xi} dt + O(N^{1/k+\xi}) \ll N^{1+1/k+\xi}. \tag{2.21}$$

Now let's look at the equation:

$$x_1 = x_2^k + x_3^k + x_4^k + x_5^k.$$

Based on (2.8), we similarly obtain the following estimate of the number of natural solutions in a hypercube with a side - $N$:

$$R_3^+(N) \ll N^{1+2/k+\xi}. \tag{2.22}$$

Let's look at the equation:

$$x_1 = x_2^k + x_3^k + x_4^k.$$

Having in mind (2.11), we similarly obtain the following estimate of the number of natural solutions in a hypercube with a side - $N$:

$$R_4^+(N) \ll N^{1+3/2k+\xi}, \qquad (2.23)$$

which corresponds to assertion 2.

Now we consider a more general Diophantine equation with an explicit variable:

$$x_1 = a_2 x_2^k + ... + a_s x_s^k, \qquad (2.24)$$

where $k, a_2, ..., a_s$ are the natural numbers.

The number of natural solutions of equation (2.24) (in a hypercube with a side - $N$) is determined by the inequality:

$$a_2 x_2^k + ... + a_s x_s^k \le N. \qquad (2.25)$$

Having in mind that $a_2 \le a_2 x_2^k, ..., a_s \le a_s x_s^k$ and (2.25), the number of natural solutions of equation (2.24) does not exceed the number of natural solutions of equation (2.15). Therefore, the asymptotic upper bound for the number of natural solutions of equation (2.24) in the hypercube with a side - $N$, determined by the formulas (2.20, 2.21, 2.22, 2.23).

Now consider the number of natural solutions of the homogeneous Diophantine equation:

$$x_1^k = x_2^k + ... + x_s^k. \qquad (2.26)$$

The equation (2.26) is Fermat's equation for the value $s = 3$ that it has no natural solutions for values $k \ge 3$, and has an infinite number of natural solutions for the value $k = 2$.

It is known that equation (2.26) has an infinite number of natural solutions for values $s = k + 1$ and $k = 1, 2, 3$. Equation (2.26) has also an infinite number of natural solutions for values $s \ge 2$ and $k = 2$.

Therefore, it is interest to find an asymptotic upper bound for the number of natural solutions of the Diophantine equation (2.26) in the above cases in a hypercube with a side - $N$.

Assertion 3

The asymptotic estimate of the number of natural solutions of equation (2.26) in a hypercube with a side - $N$ (with values $s - 1 > 2^k$) is:

$$R_s^+(N) = O(N^{(s-1)/k-1}).$$

Proof

Let's put a value $x_1^k = n$ in equation (2.26), and then we get the equation:

$$x_2^k + ... + x_s^k = n. \qquad (2.27)$$

Equation (2.27) (for a specific natural value $n$) either has no solutions or has a finite number of natural solutions. Having in mind (2.14) for values $s - 1 > 2^k$, an asymptotic upper bound for the number of natural solutions of this equation will be:

$$|R^+_{s-1}(n)| << n^{\frac{s-1}{k} - 1 + \xi}. \qquad (2.28)$$

Having in mind [7], an asymptotic upper bound for the number of natural solutions of equation (2.26) in a hypercube with a side $-N$ is equal to the sum of the solutions of equation (2.27), taking into account the probability that the sum $x_2^k + ... + x_s^k$ is equal to the $k$-th power of the natural $x_1$ - $p(n)$. In this case, the normalization condition must be fulfilled:

$$\lim_{N \to \infty} \sum_{n=1}^{N} p(n) = 1. \qquad (2.29)$$

We show that in this case:

$$p(n) = \frac{n^{1/k - 1}}{(N^{1/k} - 1)k}. \qquad (2.30)$$

We verify condition (2.29) in case (2.30):

$$\sum_{n=1}^{N} p(n) = \sum_{n=1}^{N} \frac{n^{1/k-1}}{(N^{1/k} - 1)k} = \frac{1}{(N^{1/k} - 1)k} \sum_{n=1}^{N} n^{1/k-1}. \qquad (2.31)$$

We get for a decreasing function:

$$\sum_{n=1}^{N} n^{1/k-1} = \int_{t=1}^{N} t^{1/k-1} dt + C + O(N^{1/k-1}) = k(N^{1/k} - 1) + C + O(N^{1/k-1}), \qquad (2.32)$$

where $C$ is the constant.

Having in mind (2.31) and (2.32) we obtain:

$$\sum_{n=1}^{N} p(n) = \frac{1}{(N^{1/k} - 1)k} \sum_{n=1}^{N} n^{1/k-1} = \frac{1}{(N^{1/k} - 1)k}(k(N^{1/k} - 1) + C + O(N^{1/k-1})) = 1 + \frac{C}{(N^{1/k} - 1)k} + O(\frac{N^{1/k-1}}{N^{1/k} - 1}) \qquad (2.33)$$

Based on (2.33), condition (2.29) is satisfied.

Having in mind (2.30), we obtain an asymptotic upper bound for the number of natural solutions of the Diophantine equation (2.26) in a hypercube with side - $N$:

$$R_s^+(N) = \sum_{n=1}^{N} \frac{n^{1/k-1} n^{(s-1)/k-1}}{(N^{1/k}-1)k} = \frac{1}{(N^{1/k}-1)k} \sum_{n=1}^{N} n^{s/k-2}. \tag{2.33}$$

We calculate the estimate (2.33):

$$R_s^+(N) = \frac{1}{(N^{1/k}-1)k} \int_1^N t^{s/k-2} dt + O(N^{s/k-2}) = O(N^{(s-1)/k-1}) \tag{2.34}$$

at $s > 2^k + 1$, which corresponds to assertion 3.

Formula (2.34) is valid for values $k = 2$ and $s \geq 6$. Equation (2.26) (for values $k = 2$ and $s = 6$) has the form:

$$x_1^2 = x_2^2 + ... + x_6^2. \tag{2.35}$$

Based on the above, equation (2.35) has an infinite number of natural solutions.

The asymptotic estimate of the number of natural solutions of equation (2.35) in a hypercube with a side - $N$ (based on (2.34)) is equal to:

$$R_6^+(N) = O(N^{3/2}). \tag{2.36}$$

Now let us consider Fermat equation:

$$x_1^k = x_2^k + x_3^k \tag{2.37}$$

at $k = 2$.

Assertion 4

The asymptotic estimate of the number of natural solutions of equation (2.37) in a hypercube with a side - $N$ is equal to:

$$R_{3,2}^+(N) = O(N^{1+\xi}).$$

Proof

Let's put $x_1^k = n$ in equation (2.37), and then we get the equation:

$$x_2^k + x_3^k = n. \tag{2.38}$$

Equation (2.38) either has no solutions or has a finite number of natural solutions.

Having in mind (2.6), an asymptotic upper bound for the number of natural solutions of equation (2.38) will be:

$$|R_2^+(n)| << n^{1/k+\xi}, \tag{2.39}$$

where $\xi$ is a small positive real number.

Following [7], an asymptotic upper bound for the number of natural solutions of equation (2.37) (in a side hypercube - $N$) is equal to the sum of the solutions of equation (2.38), taking into account the probability that the sum $x_2^k + x_3^k$ is equal to the $k$-degree of the natural $x_1$ - $p(n)$. The value for probability $p(n)$ is determined by the formula (2.30). Therefore, having in mind (2.30) and (2.39), we obtain an asymptotic upper bound for the number of natural solutions of the Diophantine equation (2.37) in a hypercube with side $N$:

$$|R_3^+(N)| \leq C \sum_{n=1}^{N} \frac{n^{1/k-1} n^{1/k+\xi}}{(N^{1/k}-1)\mathrm{k}} = \frac{C}{(N^{1/k}-1)\mathrm{k}} \sum_{n=1}^{N} n^{2/k-1+\xi}, \tag{2.40}$$

where $C$ is the constant.

We calculate the estimate (2.40):

$$|R_3^+(N)| \leq \frac{C}{(N^{1/k}-1)\mathrm{k}} \int_1^N t^{2/k-1+\xi} dt + O(\mathrm{N}^{2/k-1+\xi}) = O(\mathrm{N}^{2/k+\xi}). \tag{2.41}$$

It is known that, with value $k = 2$, the Fermat equation has an infinite number of natural solutions.

The asymptotic estimate of the number of natural solutions of the equation $x_1^2 = x_2^2 + x_3^2$ in a hypercube with a side $N$ (based on (2.41)) is equal to:

$$R_{3,2}^+(N) = O(N^{1+\xi}), \tag{2.42}$$

which corresponds to assertion 4.

Let's consider the equation:

$$x_1^k = x_2^k + x_3^k + x_4^k \tag{2.43}$$

for $k = 2, 3$.

Assertion 5

The asymptotic estimate of the number of natural solutions of equation (2.43) for values $k = 2, 3$ in a hypercube with the side - $N$ is:

$$R_4^+(N) = O(N^{5/2k+\xi}).$$

Proof

The corresponding equation:

$$x_2^k + x_3^k + x_4^k = n \tag{2.44}$$

based on (2.11) has the following upper bound for the number of natural solutions:

$$|R_3(n)| \ll n^{3/2k+\xi}. \tag{2.45}$$

Similarly, having in mind (2.45), we obtain an asymptotic upper bound for the number of natural solutions of equation (2.43) in a hypercube with side - $N$:

$$|R_4^+(N)| \leq C \sum_{n=1}^{N} \frac{n^{1/k-1} n^{3/2k+\xi}}{(N^{1/k}-1)k} = \frac{C}{(N^{1/k}-1)k} \sum_{n=1}^{N} n^{5/2k-1+\xi}. \tag{2.46}$$

We calculate the estimate (2.46):

$$|R_4^+(N)| \leq \frac{C}{(N^{1/k}-1)k} \int_1^N t^{5/2k-1+\xi} dt + O(N^{5/2k-1+\xi}) = O(N^{5/2k+\xi}), \tag{2.47}$$

which corresponds assertion 5.

.

Based on the above, equation (2.43) has an infinite number of solutions for values $k = 2, 3$.

The asymptotic estimate of the number of natural solutions of the equation $x_1^2 = x_2^2 + x_3^2 + x_4^2$ in a hypercube with a side - $N$ (based on (2.47)) is equal to:

$$R_{4,2}^+(N) = O(N^{5/4+\xi}). \qquad (2.48)$$

The asymptotic estimate of the number of natural solutions of the equation $x_1^3 = x_2^3 + x_3^3 + x_4^3$ in a hypercube with a side - $N$ (based on (2.47)) is equal to:

$$R_{4,3}^+(N) = O(N^{5/6+\xi}). \qquad (2.49)$$

It remains to consider the equation:

$$x_1^k = x_2^k + x_3^k + x_4^k + x_5^k \qquad (2.51)$$

for value $k = 2$.

Assertion 6

The asymptotic estimate of the number of natural solutions of the equation $x_1^k = x_2^k + x_3^k + x_4^k + x_5^k$ with value $k = 2$ in the hypercube with a side - $N$ is:

$$R_{5,2}^+(N) = O(N^{3/2+\xi}).$$

Proof

Having in mind (2.8) the corresponding equation:

$$x_2^k + x_3^k + x_4^k + x_5^k = n \qquad (2.51)$$

has the following upper bound for the number of natural solutions:

$$|R_4(n)| << n^{2/k+\xi}. \qquad (2.52)$$

Similarly, based on (2.52), we obtain an asymptotic upper bound for the number of natural solutions of equation (2.50) in a hypercube with side - $N$:

$$|R_5^+(N)| \leq C \sum_{n=1}^{N} \frac{n^{1/k-1} n^{2/k+\xi}}{(N^{1/k}-1)k} = \frac{C}{(N^{1/k}-1)k} \sum_{n=1}^{N} n^{3/k-1+\xi}. \qquad (2.53)$$

We calculate the estimate for (2.53):

$$|R_5^+(N)| \leq \frac{C}{(N^{1/k}-1)k} \int_1^N t^{3/k-1+\xi} dt + O(N^{3/k-1+\xi}) = O(N^{3/k+\xi}). \qquad (2.54)$$

Based on the above, equation (2.50) has an infinite number of solutions for value $k = 2$.

Having in mind (2.54) the asymptotic estimate of the number of natural solutions of the equation $x_1^2 = x_2^2 + x_3^2 + x_4^2 + x_5^2$ in a hypercube with a side - $N$ is equal to:

$$R_{5,2}^+(N) = O(N^{3/2+\xi}), \qquad (2.55)$$

which corresponds assertion 6.

Now consider the number of natural solutions of the inhomogeneous diagonal equation:

$$x_1 = \sum_{j=2}^{s} x_j^{k_j}, \qquad (2.56)$$

where $2 \leq k_1 \leq k_2 \leq ... \leq k_s$ and $n$ are natural numbers, $x_j$ are natural numbers and $\sum_{j=1}^{s} k_j^{-1} \geq 1$.

First, let us consider the equation:

$$\sum_{j=1}^{s} x_j^{k_j} = n. \qquad (2.57)$$

Not all natural values are representable in the form (2.57), however, almost all are representable in this form for the cases [7]:

$$k_1 = 2, k_2 = k_3 = 3; k_1 = 2, k_2 = 3, k_3 = 4; k_1 = 2, k_2 = 3, k_3 = 5.$$

Assertion 7

The asymptotic estimate of the number of natural solutions of equation (2.57) for the indicated cases in a hypercube with a side - $N$ is:

$$R_s^+(N) = O(N^{\sum_{j=2}^{s} 1/k_j + \xi}).$$

Proof

Following Vaughn and Vinogradov, the number of natural solutions satisfying equation (2.57) is equal to:

$$R_s^+(n) = \int_0^1 \prod_{i=1}^s f_i(x) e^{-2\pi x n} dx, \qquad (2.58)$$

where $f_i(x) = \sum_{m=1}^{[n^{1/k_i}]} e^{2\pi x m^{k_i}}$.

Using the Hardy – Littlewood method, it was shown in [7] that for [2.58] in the indicated cases, the following asymptotic upper bound holds:

$$R_s^+(n) << n^{\sum_{j=1}^s 1/k_j - 1 + \xi}, \qquad (2.59)$$

where $\xi$ is a small positive number.

Denote the number of natural solutions of the equation:

$$\sum_{j=2}^s x_j^{k_j} = n \qquad (2.60)$$

as $R_{s-1}^+(n)$.

Then, using (2.59), (2.60), we find the estimate of the number of natural solutions of equation (2.57) in the hypercube with side - $N$:

$$R_s^+(N) = \sum_{n=1}^N R_{s-1}^+(n) << \sum_{n=1}^N n^{\sum_{j=2}^s 1/k_j - 1 + \xi}. \qquad (2.61)$$

Having in mind that the function $n^{\sum_{j=2}^s 1/k_j - 1 + \xi}$ is not decreasing, we obtain (based on (2.61)):

$$R_s^+(N) \leq C \sum_{n=1}^N n^{\sum_{j=2}^s 1/k_j - 1 + \xi} = C \int_1^N t^{\sum_{j=2}^s 1/k_j - 1 + \xi} dt + O(N^{\sum_{j=2}^s 1/k_j - 1 + \xi}) = O(N^{\sum_{j=2}^s 1/k_j + \xi}), \qquad (2.62)$$

which corresponds to assertion 7.

Let us consider examples of determining an asymptotic upper bound for the number of natural solutions of the inhomogeneous diagonal equation (2.56) in hypercube with the side - $N$.

Let's consider the equation:

$$x_1 = x_2^2 + x_3^3 + x_4^5. \tag{2.63}$$

Having in mind (2.62), the following asymptotic upper bound for the number of natural solutions of equation (2.63) in a hypercube with a side - $N$ holds:

$$R_4^+(N) = O(N^{1/2+1/3+1/5+\xi}) = O(N^{31/30+\xi}).$$

It is possible the case, when the value (in equation (2.56)) is:

$$\sum_{j=2}^{s} 1/k_j = 1. \tag{2.64}$$

The number of natural solutions of the equation $R_3^+(n)$ is considered in [7]:

$$x_1^2 + x_2^3 + x_3^6 = n, \tag{2.65}$$

which corresponds to the case (2.64) and states that the asymptotic upper bound for the number of natural solutions of the equation:

$$x_1^2 + x_2^3 + x_3^6 = x_4$$

in the hypercube with side - $N$:

$$R_4^+(N) = \sum_{n=1}^{N} R_3^+(n) = \Gamma(3/2)\Gamma(4/3)\Gamma(7/6)N + O(N^{5/6}), \tag{2.66}$$

where $\Gamma(x)$ is the value of the gamma function at the point $x$ and $\Gamma(3/2)\Gamma(4/3)\Gamma(7/6) = 0,73...$

Estimate (2.66) corresponds to the formula (2.62).

3. ESTIMATION OF THE NUMBER OF NATUEAL (INTEGER) SOLUTIONS OF DIOPHANTINE EQUATIONSS USUNG THEIR PARAMETRIC SOLUTIONS

Suppose that the Diophantine equation has the following parametric solutions in natural (integer) numbers:

$$x_1 = \varphi_1(t_1,...,t_s),..., x_s = \varphi_s(t_1,...,t_s).$$

To make an estimate of the number of natural solutions in a hypercube with a side - $N$, you need to take the coordinate with the maximum value, for example $x_s = \varphi_s(t_1,...,t_s)$, and solve the inequality $\varphi_s(t_1,...,t_s) \leq N$.

This is quite difficult in the general case for multi-parameter solutions.

It is quite simple for one-parameter solutions. It is selected the coordinate with the maximum value of the degree of the polynomial. Suppose, that this is a polynomial of $k$ degree, then there will be the following asymptotic upper bound for the number of natural solutions in the hypercube with the side - $N$ - $O(N^{1/k})$.

For example, the equation:

$$x_1^3 + x_2^3 + x_3^3 = 1 \tag{3.1}$$

has an infinite number of integer solutions.

One-parameter solution of equation (3.1) has the form:

$$x_1 = 9a^4, x_2 = 1 - 9a^3, x_3 = 3a - 9a^4. \tag{3.2}$$

The maximum degree of the polynomial for the parametric solution (3.2) is $k = 4$. Therefore, the asymptotic upper bound for the number of natural solutions of equation (3.1) in the hypercube with a side - $N$ for solution (3.2) is equal to:

$$R_3^+(N) = O(N^{1/4}). \tag{3.3}$$

This is more difficult to make for a two-parameter solution, but in some cases it is possible.

Let's look at the Fermat equation of the second order:

$$x_1^2 + x_2^2 = x_3^2. \tag{3.4}$$

A particular two-parameter solution of equation (3.4) in natural numbers has the form:

$$x_1 = a^2 - b^2, x_2 = 2ab, x_3 = a^2 + b^2, (b > a). \tag{3.5}$$

Naturally, the coordinate $x_3 = a^2 + b^2$ takes the maximum value in solution (3.5), therefore, to estimate the number of natural solutions in a hypercube with a side - $N$, you need to find the number of natural solutions satisfying the inequalities:

$$a^2 + b^2 \leq N \text{ and } a > b. \tag{3.6}$$

Having in mind that (3.6) satisfies the number of natural points in a sector with radius $r = N^{1/2}$, then for solution (3.5) we obtain the following estimate:

$$R_3^+(N) = \frac{\pi(N^{1/2})^2}{8} = \frac{\pi N}{8} = O(N). \tag{3.7}$$

If we consider the general three-parameter solution of equation (3.4):

$$x_1 = (a^2 - b^2)c, \, x_2 = 2abc, \, x_3 = (a^2 + b^2)c, (b > a), \tag{3.8}$$

then if the third parameter $c$ changes, the values of the two-parameter solution begin to repeat. Therefore, the method of estimating the number of solutions by the volume of the figure cannot be used and the methods described above must be used in this case.

For example, an asymptotic estimate of the number of natural solutions of equation (3.4) (based on (2.42)) is:

$$R_{3,2}^+(N) = O(N^{1+\xi}),$$

those practically corresponds to the estimate of the two-parameter solution (3.7).

4. CONCLUSION AND SUGGESTIONS FOR FURTHER WORK

The next article will continue to study the asymptotic behavior of the number of integer solutions for the algebraic Diophantine equation.

5. ACKNOWLEDGEMENTS

Thanks to everyone who has contributed to the discussion of this paper.